\documentclass[11pt,a4paper]{amsart}

\usepackage[OT2,T1]{fontenc}

\usepackage[english]{babel}
\usepackage[utf8x]{inputenc}
\usepackage{amssymb,url,xspace,euscript}
\usepackage{enumerate}
\usepackage[arrow,matrix]{xy}
\xyoption{all}
\usepackage{color}
\usepackage[colorlinks,linkcolor=blue,urlcolor=blue,linktocpage]{hyperref}
\usepackage{mathtext}
\usepackage{amsmath,amscd}
\newtheorem{theo}{Theorem}[section]
 {Lemma}
{Definition}
{Proposition}
\newtheorem{coro}[theo]%
{Corollary}
{Definition-Proposition}
{Conjecture}
\newtheorem{rema}[theo]%
{Remark}
{Example}

\newtheorem{quest}[theo]%
{Question}

\newcommand{\finpreuve}{\mbox{} \hfill \mbox{$\Box$}}

\DeclareSymbolFont{cyrletters}{OT2}{wncyr}{m}{n}

\DeclareMathSymbol{\Sha}{\mathalpha}{cyrletters}{"58}

\def\F{{\rm F}}

\def\F{{\mathbb F}}
\def\N{{\mathbb N}}

\def\ds{\displaystyle}

\begin{document}

\date{}
\title[]{On the Number of Rational Points of Jacobians over Finite Fields}

\author{Philippe Lebacque}
\address{
Philippe Lebacque
\newline \indent
Laboratoire de Math\'ematiques de Besan\c con, UFR Sciences et techniques 
\newline \indent
16, route de Gray 25 030 Besan\c con, France
}
\email{philippe.lebacque@univ-fcomte.fr}

\author{Alexey Zykin}
\address{
Alexey Zykin
\newline \indent
Laboratoire GAATI, Universit\'e de la Polyn\'esie fran\c caise, BP 6570 --- 98702 Faa'a, Tahiti, Polyn\'esie fran\c caise
\newline \indent
Department of Mathematics of the National Research University Higher School of Economics
\newline \indent
AG Laboratory NRU HSE 
\newline \indent
Laboratoire Poncelet (UMI 2615)
\newline \indent
Institute for Information Transmission Problems of the Russian Academy of Sciences
}
\email{alzykin@gmail.com}

%\keywords{Galois cohomology, etale cohomology, restricted ramification}
%\subjclass[2010]{11R34, 11R37}
\thanks{The two authors were partially supported by ANR Globes ANR-12-JS01-0007-01 and the second author by AG Laboratory NRU HSE, RF government grant, ag.  11.G34.31.0023.}

\begin{abstract}{In this article we prove lower and upper bounds for class numbers of algebraic curves defined over finite fields. These bounds turn out to be better than most of the previously known bounds obtained using combinatorics. The methods used in the proof are essentially those from the explicit asymptotic theory of global fields. We thus provide a concrete application of effective results from the asymptotic theory of global fields and their zeta-functions.}
\end{abstract}

\maketitle

\section{Introduction}

\subsection{Notation}
We introduce the following notation. Let \\
\begin{tabular}{ll}
$X$ & be a smooth projective absolutely irreducible curve defined over $\F_{q}$,\\
$g$ & the genus of $X,$\\
$K$& the function field of $X,$\\
$\Phi_{{q^f}}$ &  the number of places of $K$ of degree $f,$\\
$h$ & the class number of $X$ (the number of $\F_q$-points of the Jacobian of $X$),\\
$Z_X(T)$ & the zeta function of $X$ which is a rational function of $T,$\\
$\omega_{i}\sqrt{q}$& the inverse roots of the numerator of the zeta function $Z_X(T),$\\
$\kappa$& the residue of $Z_X(q^{-s})=\zeta_X(s)$ at $s=1,$\\
$\log$& the Neperian logarithm $\log_e$.
\end{tabular}

\flushleft By a \textit{curve} we always mean a smooth projective absolutely irreducible curve.
%{\color{red}(essayer de le changer en log en base $q$! D'ailleurs, $h$ est la dimension d'un espace vectoriel sur $\F_q$ donc on a tout intérêt! Le problème est qu'il faut faire gaffe aux séries, je m'en occupe plus tard.}
\subsection{Existing lower bounds for the class number}

Our goal is to provide estimates for the number of rational points on the Jacobian of a smooth projective curve that use the information on the number of points on this curve defined over $\F_q$ or over its extensions. The starting point for all such estimates is the interpretation of the class number as the value at $1$ of the numerator of the zeta function of the curve. In order to estimate it, one uses properties of the zeta function, such as its functional equation and the Riemann Hypothesis (Weil bounds). 

From the work of Weil, we know that the class number $h$ of a smooth projective absolutely irreducible curve $X$ of genus $g$ defined over $\F_q$ is bounded by :
$$\left(\sqrt{q}-1\right)^{2g}\leq h\leq \left(\sqrt{q}+1\right)^{2g}.$$ 

Considerable efforts have been devoted to sharpening these bounds. Let us quote some works in this direction. Lachaud and Martin-Deschamps first obtained the lower bound
$$h\geq h_{LMD}=q^{g-1}\frac{(q-1)^2}{(q+1)(g+1)},$$ 
using a formula which is a consequence of the functional equation for the zeta function: 
$$\ds h=\frac{\sum\limits_{n=0}^{g-1}A_n+\sum\limits_{n=0}^{g-2}q^{g-1-n}A_n}{\sum\limits_{i=1}^g|1-\omega_i\sqrt{q}|^2},$$ where $A_n$ is the number of effective divisors of degree $n$ on $X.$ Ever since, methods from combinatorics were used to give good bounds for the numerator and the denominator of this fraction.

In a series of papers \cite{BR},\cite{BRT}, Ballet, Rolland, and Tutdere used this approach in order to prove rather elaborate lower bounds on $h$. Some of these bounds turn out to be asymptotically optimal when $g\to \infty,$ meaning that they converge to the lower bound from the generalized Brauer--Siegel theorem for function fields (\cite{TVN}). The best of their lower bounds is given by the following theorem:

\begin{theo}[Ballet--Rolland--Tutdere]
\label{BRTBound}
Let $X/\F_q$ be a curve defined over $\F_q$ of genus $g\geq 2$ and of class number $h.$ Let $D_1$, $D_2$ be finite sets of integers, $(\ell_r)_{r\in D_1},$ $(m_r)_{r\in D_2}$ be families of integers such that
\begin{enumerate}
\item $D_1\subseteq\{1,\dots,g-1\};$
\item $D_2\subseteq\{1,\dots,g-2\};$
\item for any $r\in D_1,$ $\Phi_{q^r}\geq 1;$ 
\item for any $r\in D_2,$ $\Phi_{q^r}\geq 1;$
\item $l_r\geq 0 \text{ and } \sum_{r\in D_1}r\ell_r\leq g-1;$
\item $m_r\geq 0 \text{ and } \sum_{r\in D_2}r m_r\leq g-2.$
\end{enumerate}
Then $h\geq h_{BRT}$ with
\begin{align*}h_{BRT}=&\frac{(q-1)^2}{(g+1)(q+1)-\Phi_{q}}\Bigg(\prod_{r\in D_1}\binom{\Phi_{q^r}+\ell_r}{\ell_r}\\ &\quad +q^g\prod_{r\in D_2}\Bigg[\left(\frac{q^r}{q^r-1}\right)^{\phi_{q^r}}
 -\Phi_{q^r}\binom{\Phi_{q^r}+m_r}{m_r}\int_{0}^{q^{-r}}\frac{(q^{-r}-t)^{m_r}}{(1-t)^{\Phi_{q^r}+m_r+1}}\,\mathrm{dt}\Bigg]\Bigg).
\end{align*}
\end{theo}

From now on we denote by $h_{BRT}$ the best possible lower bound from this theorem.

In a recent article dealing with estimates for the number of points on general abelian varieties, Aubry, Haloui, and Lauchaud \cite{AHL} obtained certain lower bounds on class numbers that can be very sharp when the curve in question has many rational points compared to its genus. However, these bounds are all rather poor from the asymptotic point of view when $g\to \infty$. Let us recall their results concerning Jacobian of curves.

\begin{theo}[Aubry--Haloui--Lachaud] Under the same assumptions as in theorem \ref{BRTBound}, we have:
\begin{enumerate}
\item $\ds h\geq M(q)^g\left(q+1+\frac{\Phi_q-(q+1)}{g}\right)^g,$ with $\ds M(q)=\frac{e\log x^{\frac{1}{x}-1}}{x^\frac{1}{x}-1},$ where \\$\ds x=\left(\frac{\sqrt{q}+1}{\sqrt{q}-1}\right)^2.$
\item $\ds h\geq \frac{q-1}{q^g-1}\Bigg[\binom{\Phi_q+2g-2}{2g-1}+\sum_{r=2}^{2g-1}\Phi_{q^r}\binom{\Phi_q+2g-2-r}{2g-1-i}\Bigg].$
\item If $\Phi_q\geq g(\sqrt{q}-1)+1$ then $$h\geq \binom{\Phi_q+g-1}{g}-q\binom{\Phi_q+g-3}{g-2}.$$
\item $\ds h\geq \frac{(q-1)^2}{(g+1)(q+1)-\Phi_q}\Bigg[\binom{\Phi_q+g-2}{g-2}+\sum_{r=0}^{g-1}q^{g-1-1}\binom{\Phi_q+r-1}{r}\Bigg].$
\end{enumerate}
\end{theo}

We denote by $h_{AHL}$ the best possible lower bound for $h$ given by $(1)-(4)$ of this theorem. We remark that the estimate $(3)$ can be very sharp when $g$ is small and $\Phi_q$ is large. We will come back to that in $\S 3.$

\subsection{}The aim of this paper is to show how the Mertens theorem and the explicit Brauer--Siegel theorem lead to improvements of these bounds in many cases, most notably when $g$ is large. This is done in $\S 2$ (Corollary \ref{cormain}). To do so we use the asymptotic theory of global fields and more precisely the technique of explicit formulae. The third section is devoted to numerical experiments. We compare the bounds in several examples provided by recursive asymptotically good towers of function fields. Finally, in the fourth section we discuss further research directions and open problems.
%We refer the reader to \cite{BRT} for a detailed discussion of these towers.
\vspace{0.5cm}

\textbf{Acknowledgements.} { We would like to thank St\'ephane Ballet, Julia Pieltant and Michael Tsfasman for helpful discussions.}

\section{Explicit formulae and their link to class numbers}

\subsection{Explicit formulae}
\label{sectexpl}
 Our starting point is the Mertens theorem (\cite{LMB}) for curves and its relation to the generalized Brauer--Siegel theorem. Our exposition differs slightly from \cite{LMB}: we take the opportunity to sharpen (and sometimes correct) the corresponding bounds. 

Let us recall Serre's explicit formulae from \cite{SRP}.

\begin{theo}[Explicit Formula]
\label{explicit}
For any sequence $(v_n)$ such that the radius of convergence $\rho$ of the series $\sum v_n t^n$ is strictly positive, put $\psi_{m,v}(t)=\sum_{n=1}^{+\infty}v_{mn} t^{mn},$ and $\psi_v(t)=\psi_{1,v}(t)$. Then for $t<	q^{-1}\rho$, we have the following explicit formula:
$$\sum_{f=1}^{+\infty}f\Phi_{q^f}\psi_{f,v}(t)=\psi_v(t)+\psi_v(qt)-\sum_{j=1}^{2g}\psi_v(\sqrt{q}\omega_{j}t).$$
\end{theo}

We choose $N\in\mathbb{N}$ and take $v_n=\frac{1}{n}$ if $n\leq N$ and $0$ otherwise. Applying theorem \ref{explicit} with $t=q^{-1}$, we we obtain the following identity:
$$S_0(N)=S_1(N)+S_2(N)+S_3(N),$$
where
\begin{align*}
S_0(N)= &\sum_{n=1}^{N}n^{-1} q^{-n} \sum_{m|n}m\Phi_{q^{m}}=\sum\limits_{{f=1}}^N\frac{1}{fq^f}|X(\F_{{q^f}})|,\\
S_1(N)= & \sum_{n=1}^N\frac{1}{n},\quad S_2(N)= \sum_{n=1}^N\frac{1}{nq^{n}},\quad S_3(N)= -\sum_{j=1}^{2g}\sum_{n=1}^N\frac{1}{n}(q^{-\frac{1}{2}}\omega_{j})^n.
 \end{align*}

We transform it in order to make the desired quantities appear. For any $N\geq 1,$
\begin{align*}
\underbrace{S_{0}-\sum_{f=1}^N\Phi_{{q^f}}\log\left(\frac{q^f}{q^f-1}\right)}_{\varepsilon_{0}}&+\sum_{f=1}^N\Phi_{q^f}\log\left(\frac{q^f}{q^f-1}\right)=\\
& =S_{1}+\underbrace{S_{2}-\log\frac{q}{q-1}}_{\varepsilon_{2}}+\log\frac{q}{q-1}\\
&\quad + \underbrace{S_{3}-\log\left(\kappa\log q\right)+\log\frac{q}{q-1}}_{\varepsilon_{3}}+\log\left(\kappa\log q\right)-\log\frac{q}{q-1}.
\end{align*}

In order to get bounds for $h$ we will not need estimates on $\varepsilon_{0}$ and $\varepsilon_{2}$, however they are useful for proving the Mertens theorem recalled later.

We have the following bounds for $\varepsilon_i$:
\begin{gather*}-\frac{c_{1}(q)}{N\,q^{\frac{N}{2}}}-\frac{c_{2}(q)g}{N\,q^{\frac{3N}{4}}} \leq\varepsilon_{0}\leq 0,\\
 -\frac{1}{(q-1)(N+1)q^N}\leq \varepsilon_{2}\leq 0,\quad 0\leq |\varepsilon_{3}|\leq \frac{2g}{(\sqrt{q}-1)(N+1)q^\frac{N}{2}},
\end{gather*}
with $$c_{1}(q)=\frac{2q(q+1)}{(q-1)^2}\leq 12 \ \text{ and } \  c_{2}(q)=\frac{2q}{q-1}\left(\frac{\sqrt{q}}{\sqrt{q}-1}+\frac{q^{\frac{3}{2}}}{q^{\frac{3}{2}}-1}\right)\leq 20.$$

Indeed, note the inequalities that hold for $x>1$ and $N>0$: 
$$0\leq \log\left(\frac{x}{x-1}\right)-\sum_{n=1}^N\frac{1}{nx^n}=\sum_{n=N+1}^{+\infty}\frac{1}{nx^n}\leq\frac{1}{(N+1)x^{N+1}}\sum_{n=0}^{+\infty}\frac{1}{x^n}\leq \frac{1}{(N+1)x^{N}(x-1)}.$$ This implies the bounds for $\varepsilon_2.$ The one for $\varepsilon_3$ is derived from the classical formula (\cite[Corollary 3.1.13]{TVN}):
$$\log\left(\kappa\log{q}\right)-\log\frac{q}{q-1}=\sum_{i=1}^{2g}\log\left(1-\frac{\omega_j}{\sqrt{q}}\right).$$
It gives us
\begin{align*}|\varepsilon_3|&=\left|-\sum_{j=1}^{2g}\sum_{n=1}^N\frac{1}{n}(q^{-\frac{1}{2}}\omega_{j})^n-\log\left(\kappa\log{q}\right)+\log\frac{q}{q-1}\right|\\ &=\left|\sum_{j=1}^{2g}\left(-\log\left(1-\frac{\omega_j}{\sqrt{q}}\right)-\sum_{n=1}^N\frac{1}{n}\left(\frac{\omega_{j}}{\sqrt{q}}\right)^n\right)\right|\\
\text{ and, since }|\omega_j|=1,\\ 
|\varepsilon_3|&\leq \sum_{j=1}^{2g}\frac{1}{(N+1)\sqrt{q}^{N}|\sqrt{q}-\omega_j|}\leq \frac{2g}{(\sqrt{q}-1)(N+1)q^\frac{N}{2}}.\end{align*}

\subsection{Bounds for the class number}

Using the calculations from the previous section and applying the class number formula 
$$\kappa\log q=\frac{h\,q^{1-g}}{q-1},$$
we get the following theorem.
\begin{theo}
\label{main}
Let $X$ be a smooth projective absolutely irreducible curve defined over $\F_q$ of class number $h.$ Then $h$ is given by the following formula valid for any $N\geq 1:$
$$\log h= g\log q+ \sum_{{f=1}}^N\frac{1}{fq^f}|X(\F_{{q^f}})|-\sum_{n=1}^N\frac{1+q^{-n}}{n}-\varepsilon_{3}(N),$$ or equivalently,
$$\log h= g\log q+ \sum_{{r=1}}^N\left(\Phi_{q^r}\sum_{f=1}^{\left\lfloor{\frac{N}{r}}\right\rfloor}\frac{1}{fq^{rf}}\right)-\sum_{n=1}^N\frac{1+q^{-n}}{n}-\varepsilon_{3}(N),$$
where
$\varepsilon_{3}(N)$ satisfies $\displaystyle|\varepsilon_{3}(N)|\leq \frac{2g}{(\sqrt{q}-1)(N+1)q^\frac{N}{2}}.$
\end{theo}
\begin{coro}(Bounds for the class number) 
\label{cormain}
The number of rational points $h$ on the Jacobian of $X$ satisfies $h_{min}(N)\leq h \leq h_{max}(N),$ where
$$h_{min}(N)=q^g\exp\left(\sum_{{f=1}}^N\frac{1}{fq^f}|X(\F_{{q^f}})|-\sum_{n=1}^N\frac{1+q^{-n}}{n}-\frac{2g}{(\sqrt{q}-1)(N+1)q^\frac{N}{2}}\right),$$
 $$h_{max}(N)= q^g\exp\left(\sum_{{f=1}}^N\frac{1}{fq^f}|X(\F_{{q^f}})|-\sum_{n=1}^N\frac{1+q^{-n}}{n}+\frac{2g}{(\sqrt{q}-1)(N+1)q^\frac{N}{2}}\right).$$
\end{coro}

\begin{rema}
The knowledge of a given (small) number of $\Phi_{q^f}$'s allows us, nevertheless, to apply corollary \ref{cormain} for any $N.$ For example, in the case of lower bounds, one can bound from below the unknown $\Phi_{q^f}$ by $0$ or by the quantities arising from the Weil bounds, depending which one is better. We thus get a family of bounds parametrized by $N$ and we can chose the best one. 
\end{rema}

\subsection{Mertens theorem and class numbers}

Putting together estimates from section \ref{sectexpl}, we find once again:
\begin{theo}(Mertens theorem, \cite{LMB})
\label{mert}
Let $X$ be a smooth projective absolutely irreducible curve of genus $g$ defined over $\F_q.$ Then
$$\sum_{f=1}^N\Phi_{q^f}\log\left(\frac{q^f}{q^f-1}\right)=\log\left(\kappa\log{q}\right)-\varepsilon_0+\varepsilon_{2}+\varepsilon_{3}-\sum_{n=1}^N\frac{1}{n}.$$
\end{theo}

For any $N\geq 1,$ we can deduce from it a weaker form of our bound, which might be easier compared to Ballet--Rolland--Tutdere's bound:
$$\log h= g\log q+\left[\sum_{f=1}^N\Phi_{q^f}\log\left(\frac{q^f}{q^f-1}\right)\right] -\sum_{n=1}^N\frac{1+q^{-n}}{n}+ \varepsilon_{0}-\varepsilon_{3}.$$

\begin{rema}
Theorem \ref{mert} implies that our bounds on $h$ are asymptotically optimal. More precisely, recall that a family of curves $\{X_i\}$ over $\F_q$ of genus $g_i\to\infty$ is asymptotically exact if the limits $\phi_{q^r}=\lim\limits_{i\to \infty} \frac{\Phi_{q^r}(X_i)}{g_i}$ exist for all $r.$ For asymptotically exact families of curves the generalized Brauer--Siegel theorem (\cite{TVN}) affirms that 
$$\lim_{i\to \infty} \frac{\log h(X_i)}{g_i} = \log q + \sum_{r=1}^{\infty} \phi_{q^r} \log \left(\frac{q^r}{q^r-1}\right).$$
We see that when $N\to \infty$ the bounds $h_{min}(N)$ and $h_{max}(N)$ from corollary \ref{cormain} converge to the right hand side of the above equality being divided by $g_i.$
\end{rema}

%{\color{red} Comparer asymptotiquement}

\section{Numerical computations}

In this section, we compare the lower bound $h_{min}(N)$ given by our theorem \ref{main} with $h_{BRT}$ and $h_{AHL}$ in the situation of recursive towers. We denote by $h_{LZ}$ the bound from theorem \ref{main} for the optimal choice of $N.$ Such a number $N$ is found by computer-aided calculations where the missing information on the number of points on a curve $X$ over $\F_{q^r}$ is obtained either from the inequality $X(\F_{q^r})\geq X(\F_{q^d})$ when $d \mid r,$ or from the Serre's bound $X(\F_{q^r})\geq q^r+1-g\lfloor 2q^{r/2}\rfloor,$ depending on which one is more precise. We follow closely \cite[Section 5]{BRT}.

Recall that a \textit{tower} of function fields over $\F_q$ is an infinite sequence $\{F_k/\F_q\}_{k\in \N}$ of function fields such that for all $k$ the ground field $\F_q$ is algebraically closed in $F_k,$ $F_k\subset F_{k+1},$ and the genus $g(F_k)\to \infty.$ A \textit{recursive tower} is a tower $\{F_{k}\}$ of function fields over $\F_{q}$ such that $F_{0}=\F_{q}(x_{0})$ is a rational function field and $F_{k+1}=F_{k}(x_{k+1})$ where $x_{k+1}$ satisfy the equation $f(x_{k},x_{k+1})=0$ for a given polynomial $f(X,Y)\in \F_{q}[X,Y].$

\subsection{The first tower of Garcia--Stichtenoth} 
\label{GarSt}
Assume that $q^r$ is a square and consider the tower $\{H_{k}\}=\mathcal{H}/\F_{q^r}$ defined recursively by the polynomial $$f(X,Y)=Y^{q^{\frac{r}{2}}}X^{q^{\frac{r}{2}-1}}+Y-X^{q^{\frac{r}{2}}}\in \F_{q}[X,Y].$$ We also consider the recursive tower $\{F_{k}\}=\mathcal{F}/\F_q$ of function fields defined by the same polynomial starting with the rational function field $\F_{q}(x_0).$ The base change of $F_{k}$ to $\F_{q^r}$ gives $H_{k}.$

We compare the numerical estimates from \cite[Section 5.1]{BRT} with what we obtain using our bound $h_{LZ}$. We take $q=2, r=2$ and consider the fields $H_2, H_3,$ and  $H_4.$ Note an error in \cite[Section 5.1]{BRT} where for $k=3$ the genus is erroneously taken to be equal to $14$ instead of $13$ (this was pointed out by Julia Pieltant). 
\medskip

\begin{center}
\begin{tabular}{|c|c|c|c|c|c|c|}
\hline
Step k & $g(H_k)$ &  $B_1(H_k)$ & $h_{BRT}$ & $h_{AHL}$ & $h_{LZ}$ & $N$ \\
\hline
  2 &    5 &     16 &    7434  & 12240 & 9230 & 10\\
 \hline
  3  &    13  &  30 &  16 911 279 581  &  16 271 525 520 & 26 274 427 880 & 33\\
 \hline
  4   &  33 &    56    & $1.43 \times 10^{25}$   &  $0.075 \times 10^{25}$ & $4.149 \times 10^{25}$ & 83 \\
 \hline
\end{tabular}
\end{center}

\medskip

Here is a similar comparison for $q=2$ and the tower $\mathcal F.$

\medskip

\begin{center}
\begin{tabular}{|c|c|c|c|c|c|c|}
\hline
Step k & $g(F_k)$ &   $B_1(F_k)$ &  $B_2(F_k)$ & $h_{BRT}$ & $h_{LZ}$ & $N$ \\
\hline
 2&   5 &     2 & 7 &  7 & 30 & 12\\
 \hline
 3  &     13 &  2 &   14  &   10453 & 42898 & 26\\
 \hline
 4   &  33 &     2   &   27  &  343 733 443 618 & 1 543 267 494 985 & 74\\
 \hline
 \end{tabular}
 \end{center}

\medskip

We notice that our bound is better than the other ones except for the case of $H_2/\F_4$ where we can not beat the $h_{AHL}$ bound. The situation changes, however, if we use some additional information on the places of $H_2/\F_4.$ Namely, one can calculate that $B_2(H_2)=0$ and $B_3(H_2)=24.$ These values give the bound $h_{LZ}=13430$ reached for $N=11.$ Using MAGMA we calculated that the exact value of the class number is $16200.$  

\subsection{The tower of Bassa--Garcia--Stichtenoth} 
\label{BaGarSt}
Consider the tower $\{H_{k}\}=\mathcal{H}/\F_{q^3}$ defined recursively by the polynomial $$f(X,Y)=(Y^q-Y)^{q-1}+1+\frac{X^{q(q-1)}}{(X^{q-1}-1)^{q-1}}\in \F_{q}[X,Y]$$ and let $\{F_{k}\}=\mathcal{F}/\F_q$ be the same recursive tower over $\F_q.$ We have the following numerical estimates for the class numbers when $q=2,$ that is over $\F_8$ for $H_k$ and over $\F_2$ for $F_k.$ The value of $h_{BRT}$ bound is taken from \cite[Section 5.1]{BRT}.

\medskip

\begin{center}
\begin{tabular}{|c|c|c|c|c|c|}
\hline
Step k & $g(H_k)$ &  $B_1(H_k)$  & $h_{BRT}$ &  $h_{LZ}$ & $N$  \\
\hline
 2  &    5 &     24  &  125 537 & 126832 & 9\\
 \hline
 3   &   13   &   48  &  $2.556 \times 10^{13}$ & $4.039 \times 10^{13}$ & 29\\
 \hline
 4    &    29 &      96  &  $2.010\times 10^{30}$ & $5.778\times 10^{30}$ & 11\\
 \hline
\end{tabular}
\end{center}

\medskip

\begin{center}
\begin{tabular}{|c|c|c|c|c|c|}
\hline
Step k & $g(F_k)$ &  $B_3(F_k)$ & $h_{BRT}$ & $h_{LZ}$ & $N$  \\
\hline
 2  &    5 &     8  &  3 & 3 & 5\\
 \hline
 3   &   13   &   16  &  771 & 1623 & 19 \\
 \hline
  4    &  29  &      32  &   212 127 395 & 751 622 136& 61\\
 \hline
\end{tabular}
\end{center}

\medskip

\subsection{Composite towers} The next example is the composite tower $\{E_k/\F_{q^2}\}$ constructed in \cite{HST}. It is obtained as a composite of the tower of Garcia and Stichtenoth from section \ref{GarSt} with a certain explicitly given function field. The details can be found in \cite[Propositions 5.11]{BRT}. The following table combines the estimates for $q^2=4.$

\medskip

\begin{center}
\begin{tabular}{|c|c|c|c|c|c|c|c|r|}
\hline
 Step k & $g(E_{k})$ &  $B_{1}(E_k)$  & $B_{2}(E_k)$ & $B_{3}(E_k)$ & $h_{BRT}$ & $h_{LZ}$ & $N$   \\
\hline
 2  &    55 &     1  &  12  &   12 &  3.657 $\times 10^{31}$ & 23.55 $\times 10^{31}$ & 14 \\
 \hline
 3   &   132   &   1  &   24  &  24 &  9.198 $\times 10^{77}$ & 121.02 $\times 10^{77}$ & 15 \\
 \hline
 \end{tabular}
\end{center}

\medskip

For another two composite towers $\{E_k/\F_2\}$ and $\{E'_k/\F_8\}$ this time based on the tower from section \ref{BaGarSt} (see \cite[Proposition 5.17]{BRT} for a detailed description) we get the following numerical data:

\medskip
\begin{center}
\begin{tabular}{|c|c|c|c|c|c|c|}
\hline
Step k & $g(E_k)$ &  $B_{3}(E_k)$  & $B_{6}(E_k)$ & $h_{BRT}$ & $h_{LZ}$ & $N$  \\
\hline
2  &    17 &     16  &  8   &    10 254 & 27563& 30\\
 \hline
3   &   49   &   32  &   16  &    1.718 $\times 10^{14}$ & 9.173 $\times 10^{14}$ & 94 \\
 \hline
\end{tabular}
\end{center}

\medskip

\begin{center}
\begin{tabular}{|c|c|c|c|c|c|c|}
\hline
Step k & $g(E'_k)$ &  $B_1(E'_k)$  & $B_2(E'_k)$ & $h_{BRT}$ & $h_{LZ}$ & $N$   \\
\hline
 2  &    17 &     48 &  24  & 1.002 $\times 10^{17}$ & 2.304 $\times 10^{17}$ & 35 \\
 \hline
 3   &   49   &   96  &   48 & 2.426  $\times 10^{48}$ & 13.08  $\times 10^{48}$ & 10 \\
 \hline
 \end{tabular}
\end{center}

\medskip

One more composite tower $E_k/\F_4$ introduced in \cite{W} (see also \cite[Proposition 5.18]{BRT}) gives us the following table:

\medskip

\begin{center}
\begin{tabular}{|c|c|c|c|c|c|c|c|}
\hline
Step k & $g(E_{k})$ &  $B_{1}(E_k)$  & $B_{2}(E_k)$ & $B_{3}(E_k)$ & $h_{BRT}$ & $h_{LZ}$ & $N$   \\
\hline
 2  &    30 &    1   &   9 &  9  &   4.625 $\times 10^{16}$ & 18.329 $\times 10^{16}$ & 52 \\
 \hline
  3   &    89  &  1   &   27  &   27&  2.236 $\times 10^{52}$ & 21.39 $\times 10^{52}$ & 16 \\
 \hline
 \end{tabular}
\end{center}

\medskip

For the composite tower $E_k/\F_9$ from \cite[Proposition 5.20]{BRT} we obtain.

\medskip

\begin{center}
\begin{tabular}{|c|c|c|c|c|c|c|}
\hline
 Step k & $g(E_k)$ &  $B_{1}(E_k)$  & $B_{2}(E_k)$ & $h_{BRT}$ & $h_{LZ}$ & $N$  \\
\hline
 2  &    15 &    36 &  4  &   8.563 $\times 10^{14}$ & 18.76 $\times 10^{14}$ & 30 \\
 \hline
 3   &   46   &   72  &   8  &   7.470 $\times 10^{45}$ & 41.64 $\times 10^{45}$ & 10 \\
 \hline
 \end{tabular}
\end{center}

\medskip

Finally, for yet another composite tower $E_k/\F_4$ from \cite[Proposition 5.22]{BRT} we get.

\medskip

\begin{center}
\begin{tabular}{|c|c|c|c|c|c|c|}
\hline
Step k & $g(E_k)$ &  $B_{1}(E_k)$  & $B_{2}(E_k)$ & $h_{BRT}$ & $h_{LZ}$ & $N$   \\
\hline
 2  &    25 &     36 &  9  &   1.415 $\times 10^{18}$ & 3.835 $\times 10^{18}$ & 56  \\
 \hline
 3   &   124   &   108  &   27  &  3.501  $\times 10^{86}$ & 36.23 $\times 10^{86}$ & 16 \\
 \hline
 \end{tabular}
\end{center}

\medskip

In all these examples with one exception we manage to improve on the previously known bounds.
%\begin{rema}
%For any $N\geq 1,$ we have also the weaker form, which might be easier compared to Ballet--Rolland--Tutdere's bound:
%$$\log h= g\log q+\left[\sum_{f=1}^N\Phi_{q^f}\log\left(\frac{q^f}{q^f-1}\right)\right] -\sum_{n=1}^N\frac{1+q^{-n}}{n}+ \varepsilon_{0}-\varepsilon_{3}.$$
%For our formulas to be precise, we need $\varepsilon_3(N)$ to be small. In order to get sharps bounds, $N$ shoud be taken at least..., which implies to have many places of degree $1,\dots,N...$
%\end{rema}

\section{Open questions}

Several natural questions arise in connection with the bounds obtained in this paper.

\begin{quest}
Is is possible to compare the bounds $h_{BRT}, h_{AHL},$ and $h_{LZ}$?
\end{quest}

We would like to have a more or less explicit description of the cases when each of the bounds is the best one. In the above examples our bound $h_{LZ}$ always turned out to be better than $h_{BRT}.$ However, we were not able to establish this fact in general. Comparing the bounds $h_{LZ}$ and $h_{BRT}$ does not seem to be straight forward in particular due to the fact that the number $N$ corresponding to optimal  $h_{min}(N)$ can vary significantly and does not correspond at all to the number of known $\Phi_{q^r}$'s.

\begin{quest}
Can one improve (or even optimize) the bound $h_{LZ}$ using different test functions in the explicit formulae?
\end{quest}

Oesterl\'e managed to get the best possible bounds for $|X(\F_{q^r})|$ available from explicit formulae using the linear programming approach (see \cite{SRP}). This technique, however, does not seem to be applicable directly in our case due to the non-linearity of the problem in question. The optimization seeming difficult, it would be interesting at least to find examples where different choice of test functions in the explicit formulae leads to better bounds than $h_{LZ}$. 

\begin{quest}
What are the analogues of the above bounds in the number field case?
\end{quest}

This question seems to be more directly accessible then the previous ones, since there are both the Mertens theorem and an explicit version of the Brauer--Siegel theorem available in the number field case (\cite{LMB}, \cite{LZ}). Nevertheless, analytic component in the proofs will certainly be more substantial and the application of the Generalized Riemann Hypothesis might be necessary in certain cases.

\bibliographystyle{plain}
\def\cprime{$'$}

\end{document}